\documentclass[11pt]{article}
\usepackage{amsmath}
\usepackage{amsfonts}
\usepackage{cite}

\numberwithin{equation}{section} \topmargin=0cm
\date{}
\def\limsup{\mathop{\overline{\rm lim}}}
\def\liminf{\mathop{\underline{\rm lim}}}

\begin{document}

\title {\bf Entire solutions of delay differential equations of Malmquist type }

\author{ Ran-Ran Zhang\ , \ Zhi-Bo Huang\footnote{Corresponding
author.}\\
\footnotesize
{Department of Mathematics, Guangdong University of Education, Guangzhou 510303, China}
\\
\footnotesize
{School of Mathematical Sciences, South China Normal University, Guangzhou 510631, China}
\\
\footnotesize {(Email: zhangranran@gdei.edu.cn\ , \ huangzhibo@scnu.edu.cn)}}

\date{}
\maketitle

{\small \noindent{\bf Abstract}\quad In this paper, we investigate
the delay differential equations of Malmquist type of form
\begin{equation*}
w(z+1)-w(z-1)+a(z)\frac{w'(z)}{w(z)}=R(z, w(z)),~~~~~~~~~~~~~~(*)
\end{equation*}
where $R(z, w(z))$ is an irreducible rational function in $w(z)$ with rational coefficients and
$a(z)$ is a rational function. We characterize all reduced forms when the equation $(*)$  admits a  transcendental entire solutions with hyper-order less than one. When we compare with the results obtained by Halburd and Korhonen[Proc.Amer.Math.Soc., forcoming],we obtain the reduced forms without the assumptions that the denominator of rational function $R(z,w(z))$ has roots that are nonzero rational functions in $z$.
The growth order and value distribution of transcendental entire solutions for the reduced forms are also investigated.
\\

\\

 }

\section{ Introduction and results }

We assume that the
reader is familiar with the
standard notations and basic results of the Nevanlinna theory, see e.g. \cite{hay2}.
Let $w$ be a meromorphic function in the complex plane. The order of growth of $w$ is denoted by $\sigma(w)$ and
the  hyper-order of $w$ is defined by
\begin{equation*}
\begin{split}
\sigma_2(w)=\limsup_{r\rightarrow\infty}\frac{\log\log T(r, w)}{\log r}.
\end{split}
\end{equation*}
For $a\in \mathbb{C}$, the deficiency in which zeros of $w-a$ are counted only once is defined by
\begin{equation*}
\Theta(a,
w)=1-\limsup_{r\rightarrow\infty}\frac{\overline{N}(r, \frac{1}{w-a})}{T(r, w)}.
\end{equation*}
Moreover, we say that a
meromorphic function $\alpha$ is a small function of $w$ if
$T(r, \alpha)=S(r, w)$, where $S(r, w)=o(T(r, w))$  as $r\rightarrow\infty$, possibly outside of an
exceptional set of finite logarithmic measure.
\vskip%
0.1in

The Malmquist type theorems concentrate upon necessary conditions for certain types of differential equations to admit a
meromorphic solution growing rapidly with respect to the coefficients. The following result is the celebrated Malmquist theorem.
\vskip%
0.1in

\noindent{\bf Theorem A} (\cite[p. 192]{lai1})\quad  \emph{Let $R(z, y)$
be rational in both arguments. If
the differential equation
\begin{equation*}
y'=R(z, y)
\end{equation*}
admits a transcendental meromorphic solution, then $y'=R(z, y)$ reduces
into a Riccati differential equation
\begin{equation*}
y'=a_0(z)+a_1(z)y+a_2(z)y^2
\end{equation*}
with rational coefficients. }
\vskip%
0.1in

Motivated by the problem of integrability of difference equations, Halburd and Korhonen \cite{HK2} obtained the following result, which indicates that
the existence of a finite-order meromorphic solution of a difference equation is a strong indicator of integrability of the equation.

\vskip%
0.1in

\noindent{\bf Theorem B} (\cite{HK2})\quad \emph{If the equation
\begin{equation}\label{1.1}
w(z+1)+w(z-1)=R(z, w(z))
\end{equation}
where $R(z, w(z))$ is rational in $w(z)$ with meromorphic coefficients in $z$, has an
admissible meromorphic solution of finite order, then either $w(z)$ satisfies a
difference Riccati equation
\begin{equation*}
w(z+1)=\frac{p(z+1)w(z)+q(z)}{w(z)+p(z)},
\end{equation*}
where $p, q \in \mathcal{S}(w)=\{f\ \mbox{meromorphic}: T(r, f)=o(T(r, w))\}$, or equation
(\ref{1.1}) can be transformed by a linear change in $w(z)$ to one of the following
equations:
\begin{equation*}
w(z+1)+w(z)+w(z-1)=\frac{\pi_1(z)z+\pi_2(z)}{w(z)}+\kappa_1(z)
\end{equation*}
\begin{equation*}
w(z+1)-w(z)+w(z-1)=\frac{\pi_1(z)z+\pi_2(z)}{w(z)}+(-1)^z\kappa_1(z)
\end{equation*}
\begin{equation*}
w(z+1)+w(z-1)=\frac{\pi_1(z)z+\pi_3(z)}{w(z)}+\pi_2(z)
\end{equation*}
\begin{equation*}
w(z+1)+w(z-1)=\frac{\pi_1(z)z+\kappa_1(z)}{w(z)}+\frac{\pi_2(z)}{w(z)^2}
\end{equation*}
\begin{equation*}
w(z+1)+w(z-1)=\frac{(\pi_1(z)z+\kappa_1(z))w(z)+\pi_2(z)}{(-1)^{-z}-w(z)^2}
\end{equation*}
\begin{equation*}
w(z+1)+w(z-1)=\frac{(\pi_1(z)z+\kappa_1(z))w(z)+\pi_2(z)}{1-w(z)^2}
\end{equation*}
\begin{equation*}
w(z+1)w(z)+w(z)w(z-1)=p(z)
\end{equation*}
\begin{equation*}
w(z+1)+w(z-1)=p(z)w(z)+q(z),
\end{equation*}
where $\pi_k(z), \kappa_k(z) \in \mathcal{S}(w)$ are arbitrary finite-order periodic functions with period $k$.
}
\vskip%
0.1in

Theorem B is a Malmquist type theorem for difference equations.
Furthermore, many other researchers (see, e.g. \cite{abl}, \cite{hei}, \cite{huang}, \cite{Rie}, \cite{zhangj},  also  \cite{zhang})
discussed the complex difference equations of Malmquist type.

Since some reductions of integrable differential-difference equations
are known to yield delay differential equations, Halburd and Korhonen \cite{hal2} discussed a
delay differential equation and obtained

\vskip%
0.1in

\noindent{\bf Theorem C} (\cite{hal2})\quad \emph{Let $w(z)$ be a non-rational meromorphic solution of }
\begin{equation}\label{1.2}
w(z+1)-w(z-1)+a(z)\frac{w'(z)}{w(z)}=R(z, w(z))=\frac{P(z, w(z))}{Q(z, w(z))},
\end{equation}
\emph{where $a(z)$ is rational, $P(z, w)$ is a polynomial in $w$ having rational coefficients in $z$,
and $Q(z, w)$ is a polynomial in $w(z)$ with roots that are nonzero rational functions of $z$
and not roots of $P(z, w)$. We call $f(z)$ is a subnormal solution of (\ref{eq1.2})if $f(z)$ is a solution of equation
(\ref{eq1.2}) with $\sigma_2(w)<1$, then}
\begin{equation*}
\deg_w(P)=\deg_w(Q)+1\leq 3\ \mbox{or} \ \deg_w(R)\leq 1.
\end{equation*}
\vskip%
0.1in

The notation $\deg_w(P)=\deg_w(P(z, w))$ is used to denote the degree of $P$ as a polynomial in $w$ and $\deg_w(R)=\max\{\deg_w(P), \deg_w(Q)\}$
is used to denote the degree of $R$ as a rational function in $w$.

In Theorem C, Halburd and Korhonen obtained necessary conditions for the equation (\ref{1.2}) to admit a non-rational meromorphic
solution  of hyper-order less than one under the assumption that ``$Q(z, w)$ has roots that are nonzero rational functions of $z$''. We pose two questions related to Theorem C.

1. Is it possible to obtain some reduction results as in Theorems A and B if the assumption that ``$Q(z, w)$ has roots that are nonzero rational functions of $z$'' of Theorem C
is dropped?

2. Is it possible to say something about the properties, including the form, the growth order and the distribution of $a$-values of solutions of the equation (\ref{1.2})?

In this paper, we answer these questions when the equation (\ref{1.2}) has a transcendental entire solution.
We first obtain the following result.

\vskip%
0.1in
\noindent{\bf Theorem 1.1}\quad  \emph{Let $R(z, w(z))\not\equiv0$ be an irreducible rational function in $w(z)$ with rational coefficients and
let $a(z)$ be a rational function. If the equation (\ref{1.2}) admits
a transcendental entire solution $w(z)$ with $\sigma_2(w)<1$, then (\ref{1.2}) reduces into
\begin{equation}\label{1.3}
w(z+1)-w(z-1)+a(z)\frac{w'(z)}{w(z)}=a_1(z)w(z)+a_0(z),
\end{equation}
or
\begin{equation}\label{1.4}
w(z+1)-w(z-1)+a(z)\frac{w'(z)}{w(z)}=\frac{a_2(z)w(z)^2+a_1(z)w(z)+a_0(z)}{w(z)},
\end{equation}
where $a_2(z), a_1(z), a_0(z)$ are rational functions.}
\vskip%
0.1in

\noindent{\bf Remark 1.1} (1) Theorem 1.1 is a reduction result which characterizes all cases where transcendental entire solutions of the equation (\ref{1.2})
of hyper-order less than one
actually may appear.
So Theorem 1.1 can be viewed as a weaker form of delay differential analogue of Malmquist theorem.

(2) Since $R(z, w(z))\not\equiv0$ is an irreducible rational function in $w(z)$, we see that at least one of $a_1(z)$ and $a_0(z)$
in (\ref{1.3}) does not vanish, and $a_0(z)$ in (\ref{1.4}) does not vanish.
\vskip%
0.1in

By Theorem 1.1, we easily get the following corollary.
\vskip%
0.1in

\noindent{\bf Corollary 1.1}\quad  \emph{Let $a(z)$ be a rational function, $P(z, w)$ be a polynomial in $w$ having rational coefficients in $z$,
and $Q(z, w)$ be a polynomial in $w(z)$ with roots that are nonzero rational functions of $z$
and not roots of $P(z, w)$.
Then the equation (\ref{1.2}) has no transcendental entire solutions with $\sigma_2(w)<1$.}
\vskip%
0.1in

Theorem C has an assumption ``$Q(z, w)$ has roots that are nonzero rational functions of $z$''
under which the equation (\ref{1.2}) has no transcendental entire solutions with $\sigma_2(w)<1$.
So Theorem 1.1 is independent of Theorem C, though Theorem C focuses on the case that (\ref{1.2}) has a non-rational meromorphic solution with $\sigma_2(w)<1$ and Theorem 1.1 focuses on the
case that (\ref{1.2}) has a transcendental entire solution with $\sigma_2(w)<1$.

\vskip%
0.1in

Examples 1 and 2 below show that the form (\ref{1.3}) in Theorem 1.1 does exist.
\vskip%
0.1in

\noindent{\bf Example 1.}\quad The equation
\begin{equation*}
w(z+1)-w(z-1)+a(z)\frac{w'(z)}{w(z)}=\frac{e(z+1)-e^{-1}(z-1)}{z}w(z)+a(z)\frac{1+z}{z}
\end{equation*}
has an entire solution $w(z)=ze^z$, where $a(z)$ is any rational function.
\vskip%
0.1in

\noindent{\bf Example 2.}\quad The equation
\begin{equation*}
w(z+1)-w(z-1)+a(z)\frac{w'(z)}{w(z)}=2\pi i a(z)
\end{equation*}
has an entire solution $w(z)=e^{2\pi iz}$, where $a(z)$ is any rational function.
\vskip%
0.1in
Examples 3--5 below show that the form (\ref{1.4}) in Theorem 1.1 does exist.

\vskip%
0.1in
\noindent{\bf Example 3.}\quad The equation
\begin{align*}
&w(z+1)-w(z-1)+a(z)\frac{w'(z)}{w(z)}\\
&=\frac{(e-e^{-1})w(z)^2+(-z(e-e^{-1})+2+a(z))w(z)+a(z)(1-z)}{w(z)}
\end{align*}
has an entire solution $w(z)=e^z+z$, where $a(z)$ is any rational function.
\vskip%
0.1in

\noindent{\bf Example 4.}\quad The equation
\begin{equation*}
w(z+1)-w(z-1)+a(z)\frac{w'(z)}{w(z)}=\frac{2\pi i a(z)w(z)-2\pi ia(z)}{w(z)}
\end{equation*}
has an entire solution $w(z)=e^{2\pi iz}+1$, where $a(z)$ is any rational function.
\vskip%
0.1in

\noindent{\bf Example 5.}\quad The equation
\begin{equation*}
w(z+1)-w(z-1)-\frac{1}{\pi i}\frac{w'(z)}{w(z)}=\frac{2z-\frac{1}{\pi i}}{w(z)}
\end{equation*}
has an entire solution $w(z)=e^{2\pi iz}+z$.
\vskip%
0.1in

Theorem 1.1 shows that in order to discuss the properties of transcendental entire solutions
of the equation (\ref{1.2}) with $\sigma_2(w)<1$, we only need to consider the equations (\ref{1.3}) and (\ref{1.4}). In this direction,
we get the following two results.

\vskip%
0.1in
\noindent{\bf Theorem 1.2}\quad  \emph{Let $a(z)$, $a_0(z)$ and $a_1(z)$ be rational functions with $a_1(z)\not\equiv0$ or $a_0(z)\not\equiv0$,
and let $w(z)$ be a transcendental entire solution
of the equation (\ref{1.3}) with $\sigma_2(w)<1$. }

\emph{(i) If $a(z)\equiv0$, then $\sigma(w)\geq 1$.}

\emph{(ii) If $a(z)\not\equiv0$, then $w(z)=H(z)e^{dz}$, where $H(z)\not\equiv0$ is a polynomial and $d\neq 0$ is some complex number. Especially, if
$a_1(z)$ is a polynomial with $a_1(z)\not\equiv\pm2i$, then $w(z)=Ce^{dz}$, where $C\in \mathbb{C}/\{0\}$; if $a_1(z)\equiv\pm2i$, then $w(z)=(C_1z+C_0)e^{(2k\pm\frac{1}{2})\pi iz}$,
where $k$ is an integer and $C_1, C_0\in \mathbb{C}$ with
$|C_1|+|C_0|\neq0$.}
\vskip%
0.1in

\noindent{\bf Theorem 1.3}\quad  \emph{Let $a(z)$, $a_2(z)$, $a_1(z)$ and $a_0(z)\not\equiv0$ be rational functions, and let $w(z)$ be a transcendental entire solution
of the equation (\ref{1.4}) with $\sigma_2(w)<1$. Then}

\emph{(i) $\sigma(w)\geq 1$;}

\emph{(ii) $\Theta(b, w)=0$ provided  that $b\in \mathbb{C}$ and $a_2(z)b^2+a_1(z)b+a_0(z)\not\equiv0$.}

\vskip%
0.1in

By Theorems 1.1--1.3 and Remark 1.1(2), we obtain the following result for the equation (\ref{1.2}).

\vskip%
0.1in

\noindent{\bf Corollary 1.2}\quad  \emph{Let $R(z, w(z))\not\equiv0$ be an irreducible rational function in $w(z)$ with rational coefficients,
let $a(z)$ be a rational function, and let $w(z)$ be a transcendental entire solution of the equation (\ref{1.2}) with $\sigma_2(w)<1$.
Then}

\emph{(i) $\sigma(w)\geq 1$.}

\emph{(ii) If $\deg_w(Q)=0$ and $a(z)\not\equiv0$, then  $w(z)$ has the form $w(z)=H(z)e^{dz}$, where $H(z)\not\equiv0$ is a polynomial and $d\neq 0$ is some complex number.}

\emph{(iii) If $\deg_w(Q)>0$, then $\Theta(0, w)=0$.}

\vskip%
0.1in

If $a_1(z)\equiv0$, then equation (\ref{1.3}) becomes
\begin{equation}\label{1.5}
w(z+1)-w(z-1)+a(z)\frac{w'(z)}{w(z)}=a_0(z),
\end{equation}
where $a(z), a_0(z)$ are rational. Quispel, Capel and Sahadevan \cite{QCS} showed that equation (\ref{1.5}) has a formal continuum limit to the first Painlev$\acute{\mbox e}$
equation
\begin{equation*}
\frac{d^2y}{dt^2}=6y^2+t,
\end{equation*}
if $a(z), a_0(z)$ are constants. Halburd and Korhonen \cite{hal2} indicated that if $a_0(z)\equiv p\pi ia(z)$, where $p\in \mathbb{N}$, then $w(z)=C\exp(p\pi i z)$, $C\neq0$, is a one-parameter family
of zero-free entire transcendental finite-order solutions of (\ref{1.5}) for any rational $a(z)$.

Thus, a natural question is: Dose the equation (\ref{1.5}) have entire solutions of infinite order?
In general, it is a difficult question to study the existence of meromorphic or entire solutions with $\sigma_2(w)\geq 1$ of equations involving shifts.
But for the special equation (\ref{1.5}), we obtain partial results about the  existence of entire solutions of infinite order.

For a meromorphic function $w$ of infinite order,  we use the notation of iterated order (see, e.g. \cite{B}) to express its rate of growth.
The iterated i-order of $w$ is defined by
\begin{equation*}
\begin{split}
\sigma_i(w)=\limsup_{r\rightarrow\infty}\frac{\log_{i} T(r, w)}{\log r},~(i=2, 3, 4, \cdots).
\end{split}
\end{equation*}
Obviously, the iterated 2-order of $w$ is the hyper-order of $w$.
\vskip%
0.1in

\noindent{\bf Theorem 1.4}\quad  \emph{Let $a(z)$ and $a_0(z)$ be rational functions with  $a(z)\not\equiv0$.
Then the equation (\ref{1.5}) has no entire solutions with finite iterated order. }

\section{ Proof of Theorem 1.1}
In order to prove our theorems, we need the following lemmas. The first of
these lemmas is a version of the difference analogue of the logarithmic derivative lemma.
\vskip%
0.1in

\noindent{\bf Lemma 2.1}(\cite{hal})\quad \emph{ Let $f(z)$ be a nonconstant meromorphic function
and $c\in\mathbb{C}$. If $\sigma_2(f)<1$ and $\varepsilon>0$, then
\begin{equation*}
m\left(r, \frac{f(z+c)}{f(z)}\right)=o\left(\frac{T(r, f)}{r^{1-\sigma_2(f)-\varepsilon}}\right)
\end{equation*}
for all $r$ outside of a set of finite logarithmic measure.}
\vskip%
0.1in

Applying logarithmic derivative lemma and Lemma 2.1 to Theorem 2.3 of \cite{lai2}, we get the following
lemma, which is a version of the difference analogue of the Clunie lemma.
\vskip%
0.1in

\noindent{\bf Lemma 2.2}\quad \emph{ Let $f(z)$ be a transcendental
meromorphic solution of hyper order $\sigma_2(f)<1$ of a differential-difference equation of the form
\begin{equation*}
U(z, f)P(z, f)=Q(z, f),
\end{equation*}
where $U(z, f)$ is a difference polynomial in $f(z)$ with small meromorphic coefficients,
$P(z, f),\ Q(z, f)$ are differential-difference polynomials in
$f(z)$ with small meromorphic coefficients,  $\deg_f(U)=n$ and
$\deg_f(Q)\leq n$.
Moreover, we assume that $U(z, f)$ contains just
one term of maximal total degree in $f(z)$ and its shifts. Then}
\begin{equation*}
m\big(r, P(z, f)\big)=S(r,f).
\end{equation*}
\vskip%
0.1in

The following lemma is a generalisation of
Borel's theorem on linear combinations of entire
functions.
\vskip%
0.1in

\noindent{\bf Lemma 2.3}(\cite[pp.69--70]{Gross} or \cite[p.82]{YangYi})\quad \emph{
Suppose that $f_1(z),f_2(z),\cdots,f_n(z)$ are meromorphic functions
and that $g_1(z),g_2(z),\cdots,g_n(z)$ are entire functions
satisfying the following conditions.}

(1) \emph{$\sum\limits_{j=1}^{n}f_j(z)e^{g_j(z)}\equiv0$;}

(2) \emph{$g_j(z)-g_k(z)$ are not constants for $1\leq j<k\leq n$;}

(3) \emph{for $1\leq j\leq n, 1\leq h<k\leq n$,}
\begin{equation*}
T(r,f_j)=o\{T(r,e^{g_h-g_k})\} \quad( r\rightarrow\infty,\ r\not\in E),
\end{equation*}
where $E\subset (1, \infty)$ is of finite linear measure or finite logarithmic measure.

\emph{Then $f_j(z)\equiv0$ $(j=1,2,\cdots,n)$.}
\vskip%
0.1in

\noindent{\bf Lemma 2.4} (\cite[p.29]{lai1})\quad  \emph{ Let $f$ be a meromorphic
function. Then for all irreducible rational functions in $f$
\begin{equation*}
R(z, f)=\frac{\sum_{i=0}^pa_i(z)f^i}{\sum_{j=0}^qb_j(z)f^j}
\end{equation*}
with meromorphic coefficients $a_i(z)$, $b_j(z)$ such that
\begin{equation*}
\left\{
\begin{array}{lll}
T(r, a_i)=S(r, f),\quad i=0,\cdots, p\\
T(r, b_j)=S(r, f),\quad j=0,\cdots, q,
\end{array}\right.
\end{equation*}
the characteristic function of $R(z, f)$ satisfies}
\begin{equation*}
T(r, R(z, f))=\max\{p, q\}T(r, f)+S(r, f).
\end{equation*}
\vskip%
0.1in

Lemma 2.4, due to Valiron and Mohon'ko, is of essential importance in the theory of complex differential, difference and differential-difference equations.
Next we prove the following lemma related to the equation (\ref{1.3}).
\vskip%
0.1in
\noindent{\bf Lemma 2.5} \emph{ Let $R(z, w(z))\not\equiv0$ be an irreducible rational function in $w(z)$ with rational coefficients such that $\deg_w(R)\leq2$,
let $a(z)\not\equiv0$ be a rational function and let $w(z)$ be a transcendental entire solution of the equation (\ref{1.2}). If $\sigma_2(w)<1$ and $w(z)$ has
finitely many zeros, then (\ref{1.2}) is of the form (\ref{1.3}),
where $a_1(z), a_0(z)$ are rational functions with $a_1(z)\not\equiv0$ or $a_0(z)\not\equiv0$.}
\vskip%
0.1in

{\bf  Proof of  Lemma 2.5}\quad By the hypotheses and Hadamard factorization theorem, we see that $w(z)$ takes the form
\begin{equation}\label{2.1}
w(z)=H(z)e^{g(z)},
\end{equation}
where $H(z)$ is a non-zero polynomial, $g(z)$ is a non-constant entire function such that $\sigma_2(w(z))=\sigma_2(e^{g(z)})=\sigma(g(z))<1$.
Substituting (\ref{2.1}) into the equation (\ref{1.2}) and setting
\begin{equation*}
s(z)=H(z+1)e^{g(z+1)-g(z)}-H(z-1)e^{g(z-1)-g(z)},
\end{equation*}
we get
\begin{equation}\label{2.2}
s(z)e^{g(z)}+a(z)\left(\frac{H'(z)}{H(z)}+g'(z)\right)=\frac{P(z, w(z))}{Q(z,w(z))}.
\end{equation}

If $s(z)\equiv0$, then by (\ref{2.2}), we get
\begin{equation*}
T\left(r, \frac{P(z, w(z))}{Q(z,w(z))}\right)=S(r, e^{g})=S(r, w).
\end{equation*}
By Lemma 2.4, we have $\deg_w(Q)=\deg_w(P)=0$. So the equation (\ref{1.2}) is of the form (\ref{1.3}), where
$a_0(z)\not\equiv0$ is a rational function and $a_1(z)\equiv0$.

If $s(z)\not\equiv0$, then we deduce from $\sigma_2(e^{g(z)})<1$ and Lemma 2.1 that
\begin{equation*}
T(r, s(z))=m(r, s(z))\leq m\left(r, \frac{e^{g(z+1)}}{e^{g(z)}}\right)+ m\left(r, \frac{e^{g(z-1)}}{e^{g(z)}}\right)+O(\log r)=S(r,e^{g}).
\end{equation*}
So by (\ref{2.2}), we get
\begin{equation*}
T\left(r, \frac{P(z, w(z))}{Q(z,w(z))}\right)\leq T(r, e^{g(z)})+S(r, e^{g})=T(r, w(z))+S(r, w).
\end{equation*}
This inequality and Lemma 2.4 give $\deg_w(P)\leq 1$ and $\deg_w(Q)\leq 1$. Thus, (\ref{2.2}) is of the form
\begin{equation}\label{2.3}
s(z)e^{g(z)}+a(z)\left(\frac{H'(z)}{H(z)}+g'(z)\right)=\frac{\widetilde{a_1}(z)H(z)e^{g(z)}+\widetilde{a_0}(z)}{\widetilde{b_1}(z)H(z)e^{g(z)}+\widetilde{b_0}(z)},
\end{equation}
where $\widetilde{a_1}(z), \widetilde{a_0}(z), \widetilde{b_1}(z), \widetilde{b_0}(z)$ are rational functions.
It follows from (\ref{2.3}) that
\begin{align*}
&s(z)\widetilde{b_1}(z)H(z)e^{2g(z)}+\left(\widetilde{b_0}(z)s(z)+\widetilde{b_1}(z)H(z)a(z)\left(\frac{H'(z)}{H(z)}+g'(z)\right)-\widetilde{a_1}(z)H(z)\right)e^{g(z)}\\
&+\widetilde{b_0}(z)a(z)\left(\frac{H'(z)}{H(z)}+g'(z))\right)-\widetilde{a_0}(z)=0.
\end{align*}
By this equality and Lemma 2.3, we obtain $\widetilde{b_1}(z)\equiv0$. So we deduce from (\ref{2.3}) that the equation (\ref{1.2}) is of the form (\ref{1.3}), where
$a_1(z)\not\equiv0$ and $a_0(z)$ are rational functions.
\vskip%
0.1in
In the final lemma, we consider the case where the equation (\ref{1.2}) reduces into (\ref{1.4}).

\vskip%
0.1in
\noindent{\bf Lemma 2.6}\emph{ Let $R(z, w(z))\not\equiv0$ be an irreducible rational function in $w(z)$ with rational coefficients such that $\deg_w(R)\leq2$,
let $a(z)\not\equiv0$ be a rational function and let $w(z)$ be a transcendental entire solution of the equation (\ref{1.2}). If $\sigma_2(w)<1$ and there exists
a rational function $r(z)\not\equiv0$ such that $w(z)+r(z)$ has
finitely many zeros, then (\ref{1.2}) is of the form (\ref{1.4}),
where $a_2(z), a_1(z), a_0(z)$ are rational functions with  $a_0(z)\not\equiv0$.}
\vskip%
0.1in

{\bf  Proof of Lemma 2.6}\quad By the hypotheses and Hadamard factorization theorem, we see that $w(z)$ takes the form
\begin{equation*}
w(z)=H(z)e^{g(z)}-r(z),
\end{equation*}
where $H(z)$ is a non-zero polynomial, $g(z)$ is a non-constant entire function such that $\sigma_2(w(z))=\sigma_2(e^{g(z)})=\sigma(g(z))<1$.
Setting
\begin{equation*}
s(z)=H(z+1)e^{g(z+1)-g(z)}-H(z-1)e^{g(z-1)-g(z)},
\end{equation*}
we have
\begin{equation}\label{2.4}
w(z+1)-w(z-1)=s(z)\frac{w(z)+r(z)}{H(z)}-r(z+1)+r(z-1),
\end{equation}
\begin{equation}\label{2.5}
w'(z)=\left(\frac{H'(z)}{H(z)}+g'(z)\right)w(z)+r(z)\left(\frac{H'(z)}{H(z)}+g'(z)-\frac{r'(z)}{r(z)}\right).
\end{equation}
Since $\deg_w(R)\leq2$, substituting (\ref{2.4}) and (\ref{2.5}) into the equation (\ref{1.2}),
we get
\begin{equation}\label{2.6}
\frac{\frac{s(z)}{H(z)}w(z)^2+t(z)w(z)+a(z)r(z)\left(\frac{H'(z)}{H(z)}+g'(z)-\frac{r'(z)}{r(z)}\right)}{w(z)}
=\frac{\widetilde{a_2}(z)w(z)^2+\widetilde{a_1}(z)w(z)+\widetilde{a_0}(z)}{\widetilde{b_2}(z)w(z)^2+\widetilde{b_1}(z)w(z)+\widetilde{b_0}(z)},
\end{equation}
where
\begin{equation*}
t(z)=\frac{s(z)r(z)}{H(z)}-r(z+1)+r(z-1)+a(z)\left(\frac{H'(z)}{H(z)}+g'(z)\right),
\end{equation*}
and $\widetilde{a_j}(z)$ and $\widetilde{b_j}(z) (j=0,1,2)$ are rational functions.
Since $\sigma_2(e^{g(z)})<1$, we deduce from Lemma 2.1 that $T(r, s(z))=S(r, w)$. So
all coefficients in (\ref{2.6}) are small functions of $w(z)$.
Since $H(z)$ is a polynomial, $r(z)$ is a rational function and $g(z)$ is a non-constant entire function, we deduce that
$g'(z)\not\equiv \frac{r'(z)}{r(z)}-\frac{H'(z)}{H(z)}$, which gives
\begin{equation}\label{2.7}
a(z)r(z)\left(\frac{H'(z)}{H(z)}+g'(z)-\frac{r'(z)}{r(z)}\right)\not\equiv0.
\end{equation}

If $s(z)\not\equiv0$, then multiplying both sides of (\ref{2.6}) by $w(z)(\widetilde{b_2}(z)w(z)^2+\widetilde{b_1}(z)w(z)+\widetilde{b_0}(z))$,
we get
\begin{align}\label{2.8}
&\widetilde{b_2}(z)\frac{s(z)}{H(z)}w(z)^4+t_3(z)w(z)^3+t_2(z)w(z)^2+t_1(z)w(z)\nonumber\\
&+\widetilde{b_0}(z)a(z)r(z)\left(\frac{H'(z)}{H(z)}+g'(z)-\frac{r'(z)}{r(z)}\right)=0,
\end{align}
where $t_j(z)(j=1,2,3)$ are all small functions of $w(z)$. By (\ref{2.7}), (\ref{2.8}) and Lemma 2.4,
we have $\widetilde{b_2}(z)\equiv0$ and $\widetilde{b_0}(z)\equiv0$. So the equation (\ref{1.2}) is of the form (\ref{1.4}), where
$a_j(z) (j=0,1,2)$ are rational functions. Recalling that $R(z, w(z))$ is irreducible, we get $a_0(z)\not\equiv0$.

If $s(z)\equiv0$, then by (\ref{2.6}) and Lemma 2.4, we have
$\deg_w(P)\leq1$ and $\deg_w(Q)\leq1$. So (\ref{2.6}) reduces into
\begin{equation}\label{2.9}
\frac{t(z)w(z)+a(z)r(z)\left(\frac{H'(z)}{H(z)}+g'(z)-\frac{r'(z)}{r(z)}\right)}{w(z)}
=\frac{\widetilde{a_1}(z)w(z)+\widetilde{a_0}(z)}{\widetilde{b_1}(z)w(z)+\widetilde{b_0}(z)}.
\end{equation}
By (\ref{2.9}) and using the same reasoning as above, we see that the equation (\ref{1.2}) is of the form (\ref{1.4}), where $a_2(z)\equiv0$
and $a_j(z) (j=0,1)$ are rational functions with $a_0(z)\not\equiv0$.
\vskip%
0.1in

{\bf Proof of Theorem 1.1}\quad First we discuss the case $a(z)\equiv0$. We deduce from (\ref{1.2}) and Lemma 2.1 that
\begin{align}\label{2.10}
T\left(r, \frac{P(z, w(z))}{Q(z,w(z))}\right)
&=T(r, w(z+1)-w(z-1))\nonumber\\
&=m(r, w(z+1)-w(z-1))\nonumber\\
&\leq m(r, w(z))+m\left(r, \frac{w(z+1)}{w(z)}\right)+m\left(r, \frac{w(z-1)}{w(z)}\right)+O(1)\nonumber\\
&=m(r, w(z))+S(r, w).
\end{align}
By (\ref{2.10}) and Lemma 2.4, we have $\deg_w(P)\leq1$ and $\deg_w(Q)\leq1$. So the equation (\ref{1.2}) has the form
\begin{equation}\label{2.11}
w(z+1)-w(z-1)=\frac{\widetilde{a_1}(z)w(z)+\widetilde{a_0}(z)}{\widetilde{b_1}(z)w(z)+\widetilde{b_0}(z)},
\end{equation}
where $\widetilde{a_j}(z)$ and $\widetilde{b_j}(z) (j=0,1)$ are rational functions.
We affirm that $\widetilde{b_1}(z)\equiv0$. Otherwise, if $\widetilde{b_1}(z)\not\equiv0$, we deduce from Lemma 2.4 and (\ref{2.11}) that
\begin{equation}\label{2.12}
T(r, w(z+1)-w(z-1))=T\left(r, \frac{\widetilde{a_1}(z)w(z)+\widetilde{a_0}(z)}{\widetilde{b_1}(z)w(z)+\widetilde{b_0}(z)}\right)=T(r, w(z))+S(r, w).
\end{equation}
Moreover, by (\ref{2.11}), we get
\begin{equation}\label{2.13}
w(z)(w(z+1)-w(z-1))=-\frac{\widetilde{b_0}(z)}{\widetilde{b_1}(z)}(w(z+1)-w(z-1))+\frac{\widetilde{a_1}(z)}{\widetilde{b_1}(z)}w(z)+\frac{\widetilde{a_0}(z)}{\widetilde{b_1}(z)}.
\end{equation}
Applying  Lemma 2.2 to (\ref{2.13}), we get
\begin{equation*}
T(r, w(z+1)-w(z-1))=m(r, w(z+1)-w(z-1))=S(r, w).
\end{equation*}
This contradicts (\ref{2.12}). So $\widetilde{b_1}(z)\equiv0$ and  the equation (\ref{1.2}) is of the form (\ref{1.3}), where $a_j(z) (j=0,1)$ are rational functions
with $a_1(z)\not\equiv0$ or $a_0(z)\not\equiv0$.

Next we discuss the case $a(z)\not\equiv0$. We deduce from (\ref{1.2}) and Lemma 2.1 that
\begin{align*}
&T\left(r, \frac{P(z, w(z))}{Q(z,w(z))}\right)\\
&=m\left(r, w(z+1)-w(z-1)+a(z)\frac{w'(z)}{w(z)}\right)+N\left(r, w(z+1)-w(z-1)+a(z)\frac{w'(z)}{w(z)}\right)\\
&\leq m(r, w(z))+N\left(r, \frac{1}{w(z)}\right)+S(r, w)\\
&\leq 2T(r, w(z))+S(r, w).
\end{align*}
By Lemma 2.4, we have $\deg_w(P)\leq2$ and $\deg_w(Q)\leq2$. So
\begin{equation}\label{2.14}
\frac{P(z, w(z))}{Q(z,w(z))}=\frac{\widetilde{a_2}(z)w(z)^2+\widetilde{a_1}(z)w(z)+\widetilde{a_0}(z)}{\widetilde{b_2}(z)w(z)^2+\widetilde{b_1}(z)w(z)+\widetilde{b_0}(z)},
\end{equation}
where $\widetilde{a_j}(z)$ and $\widetilde{b_j}(z) (j=0,1,2)$ are rational functions.

If $w(z)$ has
finitely many zeros, then Lemma 2.5 shows that the equation (\ref{1.2}) is of the form (\ref{1.3}).
If there exists
a rational function $r(z)\not\equiv0$ such that $w(z)+r(z)$ has
finitely many zeros, then lemma 2.6 shows that the equation (\ref{1.2}) is of the form (\ref{1.4}).

Now we assume that $w(z)$ has infinitely many zeros and  $w(z)+r(z)$ also has infinitely many zeros for any rational function $r(z)\not\equiv0$.
Suppose that $z_0$ is a zero of $w(z)$ and that neither $a(z)$ nor any of the coefficients in $\frac{P(z, w(z))}{Q(z,w(z))}$
has a zero or a pole at $z_0$. If $\widetilde{b_0}(z)\not\equiv0$, then $z_0$ is a simple pole of
$w(z+1)-w(z-1)+a(z)\frac{w'(z)}{w(z)}$ and a finite value of
$\frac{P(z, w(z))}{Q(z, w(z))}$, a contradiction. So $\widetilde{b_0}(z)\equiv0$.

If $\widetilde{b_2}(z)\not\equiv0$ and  $\widetilde{b_1}(z)\equiv0$, then by (\ref{1.2}) and (\ref{2.14}) we have
\begin{equation}\label{2.15}
w(z+1)-w(z-1)+a(z)\frac{w'(z)}{w(z)}=\frac{\widetilde{a_2}(z)w(z)^2+\widetilde{a_1}(z)w(z)+\widetilde{a_0}(z)}{\widetilde{b_2}(z)w(z)^2}.
\end{equation}
Since the right hand side of (\ref{2.15}) is irreducible in $w(z)$, we see that $\widetilde{a_0}(z)\not\equiv0$.
Choose a zero $z_0$ of $w(z)$ as above. Then we see that $z_0$ is a simple pole of the left hand side of (\ref{2.15}) and
a multiple pole of the right hand side of (\ref{2.15}), a contradiction.

If $\widetilde{b_2}(z)\not\equiv0$ and  $\widetilde{b_1}(z)\not\equiv0$, then by (\ref{1.2}) and (\ref{2.14}) we have
\begin{equation}\label{2.16}
w(z+1)w(z)-w(z-1)w(z)+a(z)w'(z)
=\frac{\frac{\widetilde{a_2}(z)}{\widetilde{b_2}(z)}w(z)^2+\frac{\widetilde{a_1}(z)}{\widetilde{b_2}(z)}w(z)+\frac{\widetilde{a_0}(z)}{\widetilde{b_2}(z)}}{w(z)+\frac{\widetilde{b_1}(z)}{\widetilde{b_2}(z)}}.
\end{equation}
Since the right hand side of (\ref{2.16}) is irreducible in $w(z)$, we see that
$\frac{\widetilde{a_2}(z)}{\widetilde{b_2}(z)}w(z)^2+\frac{\widetilde{a_1}(z)}{\widetilde{b_2}(z)}w(z)+\frac{\widetilde{a_0}(z)}{\widetilde{b_2}(z)}$ and
$w(z)+\frac{\widetilde{b_1}(z)}{\widetilde{b_2}(z)}$ have at most finitely many common zeros.
Furthermore,
$w(z)+\frac{\widetilde{b_1}(z)}{\widetilde{b_2}(z)}$ has infinitely many zeros. So we can choose a zero $z_1$ of $w(z)+\frac{\widetilde{b_1}(z)}{\widetilde{b_2}(z)}$ such that
neither $a(z)$ nor $\frac{\widetilde{a_2}(z)}{\widetilde{b_2}(z)}w(z)^2+\frac{\widetilde{a_1}(z)}{\widetilde{b_2}(z)}w(z)+\frac{\widetilde{a_0}(z)}{\widetilde{b_2}(z)}$
has a zero or a pole at $z_1$.
So $z_1$ is a pole of the right hand side of (\ref{2.16}) and a finite value of the left hand side of (\ref{2.16}), a contradiction.

From the above discussion, we see that $\widetilde{b_0}(z)\equiv0$ and $\widetilde{b_2}(z)\equiv0$. So the equation (\ref{1.2}) is of the form (\ref{1.4}),
where $a_j(z) (j=0,1,2)$ are rational functions with $a_0(z)\not\equiv0$.

\section{ Proof of Theorem 1.2}
In order to prove Theorem 1.2, we need the following lemmas for linear difference equations.
\vskip%
0.1in

\noindent{\bf Lemma 3.1} (\cite{chenshon})\quad\emph{ Let $P_n(z),\cdots, P_0(z)$ be polynomials such that
$P_n(z)P_0(z)\not\equiv0$ and satisfy $P_n(z)+\cdots+P_0(z)\not\equiv0$. Then every transcendental meromorphic solution $f(z)$
of the equation
\begin{equation*}
P_n(z)f(z+n)+P_{n-1}(z)f(z+n-1)+\cdots+P_0(z)f(z)=0
\end{equation*}
satisfies $\sigma(f)\geq1$.}
\vskip%
0.1in

\noindent{\bf Lemma 3.2} (\cite{chenshon})\quad\emph{ Let $F(z), P_n(z),\cdots, P_0(z)$ be polynomials such that
$F(z)P_n(z)P_0(z)\not\equiv0$. Then every transcendental meromorphic solution $f(z)$
of the equation
\begin{equation*}
P_n(z)f(z+n)+P_{n-1}(z)f(z+n-1)+\cdots+P_0(z)f(z)=F(z)
\end{equation*}
satisfies $\sigma(f)\geq1$.}
\vskip%
0.1in

{\bf Proof of Theorem 1.2}\quad
(i) Since $a(z)\equiv0$, we get from (\ref{1.3}) that
\begin{equation*}
w(z+2)-a_1(z+1)w(z+1)-w(z)=a_0(z+1).
\end{equation*}
If $a_0(z)\not\equiv0$, then by Lemma 3.2, we have $\sigma(w)\geq1$. If $a_0(z)\equiv0$, then $a_1(z)\not\equiv0$.
So by Lemma 3.1, we have $\sigma(w)\geq1$.

(ii) Since $a(z)\not\equiv0$, we see from (\ref{1.3}) that $w(z)$ has only finitely many zeros.
By Hadamard factorization theorem, $w(z)$ takes the form
\begin{equation}\label{3.1}
w(z)=H(z)e^{g(z)},
\end{equation}
where $H(z)$ is a non-zero polynomial, $g(z)$ is a non-constant entire function such that $\sigma_2(w(z))=\sigma_2(e^{g(z)})=\sigma(g(z))<1$.
Substituting (\ref{3.1}) into (\ref{1.3}), we get
\begin{equation}\label{3.2}
-a_1(z)H(z)e^{g(z)}+H(z+1)e^{g(z+1)}-H(z-1)e^{g(z-1)}=a_0(z)-a(z)\left(\frac{H'(z)}{H(z)}+g'(z)\right).
\end{equation}
Since $\sigma(g(z))=\sigma_2(e^{g(z)})<1$, we have $\liminf\limits_{r\rightarrow\infty}\frac{T(r, g)}{r}=0.$
If $g(z)$ is a transcendental entire function, we see from \cite[p. 101]{whi} that $g(z+1)-g(z)$, $g(z+1)-g(z-1)$ and $g(z-1)-g(z)$ are all transcendental entire functions.
Applying Lemma 2.3 to (\ref{3.2}), we get $H(z+1)\equiv0$, a contradiction. So $g(z)$ must be a polynomial.
If $\deg g(z)\geq2$, then $\deg(g(z+1)-g(z))=\deg(g(z+1)-g(z-1))=\deg(g(z-1)-g(z))\geq1$. Using Lemma 2.3 again, we also get
a contradiction. So  $\deg g(z)=1$ and $w(z)$ has the form
\begin{equation}\label{3.3}
w(z)=H(z)e^{dz},
\end{equation}
where $d\neq 0$ is some complex number.
Substituting (\ref{3.3}) into (\ref{1.3}), we get
\begin{equation}\label{3.4}
e^{dz}(H(z+1)e^d-H(z-1)e^{-d})=a_1(z)H(z)e^{dz}+a_0(z)-a(z)\left(\frac{H'(z)}{H(z)}+d\right).
\end{equation}
Applying Lemma 2.3 to (\ref{3.4}), we get
\begin{equation}\label{3.5}
H(z+1)e^d-H(z-1)e^{-d}=a_1(z)H(z).
\end{equation}
If $a_1(z)$ is a polynomial, we deduce from (\ref{3.5}) that $a_1(z)$ must be a constant.
Let
\begin{equation}\label{3.6}
H(z)=c_nz^n+c_{n-1}z^{n-1}+\cdots, \  (c_n\neq0).
\end{equation}
If $n\geq2$, then by (\ref{3.5}) and (\ref{3.6}), we get
\begin{equation}\label{3.7}
\left\{
\begin{array}{lll}
e^d-e^{-d}=a_1(z),\\
c_{n-1}(e^d-e^{-d})+nc_n(e^d+e^{-d})=a_1(z)c_{n-1},\\
c_{n-2}(e^d-e^{-d})+\frac{n(n-1)}{2}c_n(e^d-e^{-d})+(n-1)c_{n-1}(e^d+e^{-d})=a_1(z)c_{n-2}.
\end{array}
\right.
\end{equation}
(\ref{3.7}) gives $e^d=e^{-d}=0$. This is impossible. So $n\leq1$.

If $a_1(z)\not\equiv \pm2i$, we must have $n=0$. Otherwise, if
$n=1$, then by (\ref{3.5}) and  (\ref{3.6}), we get
\begin{equation*}
\left\{
\begin{array}{lll}
e^d-e^{-d}=a_1(z),\\
e^d+e^{-d}=0.
\end{array}
\right.
\end{equation*}
This gives $e^d=\pm i$ and $a_1(z)=\pm2i$, contradicts $a_1(z)\not\equiv \pm2i$.
So $H(z)$ must be a constant and $w(z)$ has the form $w(z)=Ce^{dz}$, where $C\neq0, C\in \mathbb{C}$.

If $a_1(z)\equiv \pm2i$, then by (\ref{3.5}) and  (\ref{3.6}), we get
\begin{equation*}
e^d-e^{-d}=\pm2i.
\end{equation*}
So $d=(2k\pm\frac{1}{2})\pi i$ and $w(z)=(C_1z+C_0)e^{(2k\pm\frac{1}{2})\pi iz}$,
where $k$ is an integer and $C_1, C_0\in \mathbb{C}$ with
$|C_1|+|C_0|\neq0$.

\section{ Proof of Theorem 1.3}
In order to prove Theorem 1.3, we need the following lemmas. The first of these lemmas is another version of difference analogue of the logarithmic derivative lemma.
\vskip%
0.1in

\noindent{\bf Lemma 4.1} (\cite{chi})\quad\emph{
Let $\eta_1$, $\eta_2$ be
two complex numbers such that $\eta_1\neq \eta_2$ and let $f(z)$ be
a finite order meromorphic function. Let $\sigma$ be the order of
$f(z)$, then for each $\varepsilon>0$, we have
\begin{equation*}
m\left(r,
\frac{f(z+\eta_1)}{f(z+\eta_2)}\right)=O(r^{\sigma-1+\varepsilon}).
\end{equation*}
}

By a careful inspection of the proof of
Theorem 2.3 in \cite{lai2} and using Lemma 4.1, we easily get the following lemma,
which is another version of the difference analogue of the Clunie lemma.
\vskip%
0.1in

\noindent{\bf Lemma 4.2}\quad \emph{ Let $f(z)$ be a transcendental
meromorphic solution of finite order $\sigma$ of a differential-difference equation of the form
\begin{equation*}
f(z)^nP(z, f)=Q(z, f),
\end{equation*}
where $P(z, f),\ Q(z, f)$ are differential-difference polynomials in
$f(z)$ with rational coefficients and
$\deg_f(Q)\leq n$. Then for each $\varepsilon>0$, we have}
\begin{equation*}
m\big(r, P(z, f)\big)=O(r^{\sigma-1+\varepsilon})+O(\log r).
\end{equation*}
\vskip%
0.1in

Applying Lemma 2.1 to Corollary 3.4
of \cite{HK1}, we easily get the following lemma,
which is a version of the difference analogue of the Mohon'ko-Mohon'ko lemma.
\vskip%
0.1in

\noindent{\bf Lemma 4.3}\quad \emph{ Let $f(z)$ be a transcendental meromorphic solution with $\sigma_2(f)<1$ of
\begin{equation*}
P(z, f)=0,
\end{equation*}
where $P(z, f)$ is a differential-difference polynomial in $f(z)$ with small meromorphic coefficients. If $P(z, a)\not\equiv0$ for a small function $a(z)$, then
\begin{equation*}
m\left(r, \frac{1}{f-a}\right)=S(r, f).
\end{equation*}}
\vskip%
0.1in

{\bf Proof of Theorem 1.3}\quad
(i) If $a(z)\equiv0$, then we deduce from (\ref{1.4}) and Lemma 2.1 that
\begin{align*}
&T\left(r, \frac{a_2(z)w(z)^2+a_1(z)w(z)+a_0(z)}{w(z)}\right)\\
&=T(r, w(z+1)-w(z-1))\\
&=m(r, w(z+1)-w(z-1))\\
&\leq m(r, w(z))+S(r, w).
\end{align*}
So by Lemma 2.4, we see that $a_2(z)\equiv0$ and (\ref{1.4}) can be written as
\begin{equation}\label{4.1}
w(z)(w(z+1)-w(z-1))=a_1(z)w(z)+a_0(z).
\end{equation}
By (\ref{4.1}), Lemma 2.2 and Lemma 2.4, we have
\begin{align*}
T(r, w(z))+S(r, w)&=T\left(r, \frac{a_1(z)w(z)+a_0(z)}{w(z)}\right)\\
&=T(r, w(z+1)-w(z-1))\\
&=m(r, w(z+1)-w(z-1))\\
&=S(r, w).
\end{align*}
This is a contradiction. So we must have $a(z)\not\equiv0$.

Assume that $\sigma(w(z))<1$. We get from (\ref{1.4}) that
\begin{equation}\label{4.2}
w(z)(w(z+1)-w(z-1)-a_2(z)w(z)-a_1(z))=-a(z)w'(z)+a_0(z).
\end{equation}
Since $w(z)$ is transcendental and $a(z)\not\equiv0$, we see that $-a(z)w'(z)+a_0(z)\not\equiv0$. So
$w(z+1)-w(z-1)-a_2(z)w(z)-a_1(z)\not\equiv0$. Applying  Lemma 4.2 to (\ref{4.2}), we get
\begin{align}\label{4.3}
&T(r, w(z+1)-w(z-1)-a_2(z)w(z)-a_1(z))\nonumber\\
&=m(r, w(z+1)-w(z-1)-a_2(z)w(z)-a_1(z))\nonumber\\
&=O(\log r).
\end{align}
Rewrite (\ref{4.2}) as
\begin{equation}\label{4.4}
-\frac{w(z+1)-w(z-1)-a_2(z)w(z)-a_1(z)}{a(z)}=\frac{w'(z)}{w(z)}-\frac{a_0(z)}{a(z)w(z)}.
\end{equation}
By (\ref{4.3}) and $a(z)$ is rational, we see that
$-\frac{w(z+1)-w(z-1)-a_2(z)w(z)-a_1(z)}{a(z)}\not\equiv0$ is rational. So
\begin{equation}\label{4.5}
-\frac{w(z+1)-w(z-1)-a_2(z)w(z)-a_1(z)}{a(z)}=Az^n(1+o(1)),
\end{equation}
where $z\rightarrow\infty$, $A\neq0$ is a constant and $n$ is an integer.
From the Wiman-Valiron theory (see  \cite[pp.28-32]{he}, \cite[p.51]{lai1} or \cite[pp.103-105]{valiron}), we have
\begin{equation}\label{4.6}
\frac{w'(z)}{w(z)}=\frac{\nu(r)}{z}(1+o(1)),
\end{equation}
where $|z|=r\not\in[0,1]\bigcup E, E\subset (1, \infty)$ is of finite logarithmic measure such that
$|w(z)|=M(r, w)$ and $\nu(r)$ denotes the central index of $w(z)$.
Substituting (\ref{4.5}) and (\ref{4.6}) into (\ref{4.4}), we have
\begin{equation}\label{4.6+1}
Az^{n+1}(1+o(1))=\nu(r)(1+o(1))-\frac{za_0(z)}{a(z)w(z)},
\end{equation}
where $|z|=r\not\in[0,1]\bigcup E$ such that
$|w(z)|=M(r, w)$.
Since $w(z)$ is transcendental, we have
\begin{equation}\label{4.6+2}
\frac{|za_0(z)|}{|a(z)|M(r,w)}\rightarrow 0,\ (r\rightarrow\infty).
\end{equation}
By (\ref{4.6+1}) and (\ref{4.6+2}), we get
\begin{equation*}
\nu(r)=|A|r^{n+1}(1+o(1)), \ (r\not\in[0,1]\bigcup E,\ r\rightarrow\infty).
\end{equation*}
So $\sigma(w(z))=\limsup\limits_{r\rightarrow\infty}\frac{\log \nu(r)}{\log r}=n+1.$ Since $\sigma(w(z))<1$ and
$n+1$ is an integer, we have $n+1\leq 0$. Thus $\nu(r)=|A|r^{n+1}(1+o(1))$ contradicts the fact that
$w(z)$ is transcendental. Therefore we proved that  $\sigma(w(z))\geq 1$.

(ii) Let
\begin{equation*}
P(z, w)=w(z)(w(z+1)-w(z-1))+a(z)w'(z)-a_2(z)w(z)^2-a_1(z)w(z)-a_0(z).
\end{equation*}
Then $P(z, b)=-a_2(z)b^2-a_1(z)b-a_0(z)\not\equiv0$. By lemma 4.3, we have
\begin{equation*}
m\left(r, \frac{1}{w(z)-b}\right)=S(r, w).
\end{equation*}
So
\begin{equation}\label{4.7}
N\left(r, \frac{1}{w(z)-b}\right)=T(r, w(z))+S(r, w).
\end{equation}

First we suppose that $b=0$. Then equation (\ref{1.4}) shows that $w(z)$ has at most finitely many multiple zeros.
So by (\ref{4.7}), we get $\Theta(0, w(z))=0$.

Next we suppose that $b\neq0$.
Let $g(z)=w(z)-b$, then $g(z)$ has infinitely many zeros.
Substituting $w(z)=g(z)+b$ into the equation (\ref{1.4}), we get
\begin{align}\label{4.8}
g(z)(g(z+1)-g(z-1))&+b(g(z+1)-g(z-1))+a(z)g'(z)\nonumber\\
&=a_2(z)g(z)^2+2ba_2(z)g(z)+a_1(z)g(z)+\psi(z),
\end{align}
where $\psi(z)=a_2(z)b^2+a_1(z)b+a_0(z)$.
Now we divide our discussion into two cases.

Case 1. $g(z+1)-g(z-1)-a_2(z)g(z)\not\equiv0$.

Rewrite (\ref{4.8}) as
\begin{align}\label{4.9}
g(z)(&g(z+1)-g(z-1)-a_2(z)g(z))\nonumber\\
&=-a(z)g'(z)-b(g(z+1)-g(z-1))+2ba_2(z)g(z)+a_1(z)g(z)+\psi(z).
\end{align}
Applying Lemma 2.2 to (\ref{4.9}), we get
\begin{equation*}
T(r, g(z+1)-g(z-1)-a_2(z)g(z))=m(r, g(z+1)-g(z-1)-a_2(z)g(z))=S(r, g).
\end{equation*}
If $g(z+1)-g(z-1)-a_2(z)g(z)\equiv \frac{\psi(z)}{b}$, then by (\ref{4.8}), we get
\begin{equation*}
g(z)(g(z+1)-g(z-1))+a(z)g'(z)=a_2(z)g(z)^2+ba_2(z)g(z)+a_1(z)g(z).
\end{equation*}
Comparing the orders of zeros of both sides of the above equality, we get a contradiction. So
$g(z+1)-g(z-1)-a_2(z)g(z)\not\equiv \frac{\psi(z)}{b}$ and
\begin{equation}\label{4.10}
N\left(r, \frac{1}{g(z+1)-g(z-1)-a_2(z)g(z)-\frac{\psi(z)}{b}}\right)=S(r, g).
\end{equation}

We denote by $N_{1}\left(r, \frac{1}{g(z)}\right)$ the counting function of those simple zeros of $g(z)$  in
$|z|<r$,
and denote by $N_{>1}\left(r, \frac{1}{g(z)}\right)$ the counting function of those multiple zeros of $g(z)$ in
$|z|<r$. If $z_0$ is a multiple zero of $g(z)$ and that none of the coefficients in (\ref{4.8}) has a zero
or a pole at $z_0$, then by (\ref{4.8}), we have
\begin{equation*}
b(g(z_0+1)-g(z_0-1))=\psi(z_0),
\end{equation*}
and so
\begin{equation}\label{4.11}
g(z_0+1)-g(z_0-1)-a_2(z_0)g(z_0)-\frac{\psi(z_0)}{b}=0.
\end{equation}
By (\ref{4.10}) and (\ref{4.11}),
we get
\begin{align*}
&N\left(r, \frac{1}{g(z)}\right)\nonumber\\
&=N_{1}\left(r, \frac{1}{g(z)}\right)+N_{>1}\left(r, \frac{1}{g(z)}\right)\nonumber\\
&\leq N_{1}\left(r, \frac{1}{g(z)}\right)+N\left(r, \frac{1}{g(z+1)-g(z-1)-a_2(z)g(z)-\frac{\psi(z)}{b}}\right)+O(\log r)\nonumber\\
&\leq \overline{N}\left(r, \frac{1}{g(z)}\right)+S(r, g),
\end{align*}
which gives $\Theta(b, w(z))=\Theta(0, g(z))=0$.

Case 2. $g(z+1)-g(z-1)-a_2(z)g(z)\equiv0$.

Substituting $b(g(z+1)-g(z-1))=ba_2(z)g(z)$ into (\ref{4.8}), we have
\begin{equation*}
g(z)( g(z+1)-g(z-1))+a(z)g'(z)=a_2(z)g(z)^2+ba_2(z)g(z)+a_1(z)g(z)+\psi(z).
\end{equation*}
Since $\psi(z)\not\equiv0$, we see from the above equality that $g(z)$ has at most finitely many multiple zeros.
So $\Theta(b, w(z))=\Theta(0, g(z))=0$.

\section{ Proof of Theorem 1.4}
In order to prove Theorem 1.4, we need the following lemma, which relates to the estimate of characteristic function of shifts of a meromorphic function.
\vskip%
0.1in

\noindent{\bf Lemma 5.1} (\cite{abl, GO}).\quad\emph{
Let $f(z)$ be a meromorphic function. For an arbitrary $c\neq0$, the
following inequalities
\begin{equation*}
\begin{split}
\big(1+o(1)\big)T\big(r-|c|, f(z)\big)\leq T\big(r, f(z+c)\big)\leq
\big(1+o(1)\big)T\big(r+|c|, f(z)\big)
\end{split}
\end{equation*}
hold as $r\rightarrow\infty$.}
\vskip%
0.1in

{\bf  Proof of Theorem 1.4.}\quad
Suppose that $w(z)$ is an entire solution of the equation (\ref{1.5}) with finite iterated
order. Then there exists
an integer $p\geq1$ such that $\sigma_p(w)=\infty$ and $\sigma_{p+1}(w)<\infty$. By
equation (\ref{1.5}), we get
\begin{equation}\label{5.1}
a_0(z)-a(z)\frac{w'(z)}{w(z)}=w(z+1)-w(z-1).
\end{equation}
By (\ref{5.1}) and Hadamard factorization theorem, we see that $w(z)$ takes the form
\begin{equation}\label{5.2}
w(z)=H(z)e^{g(z)},
\end{equation}
where $H(z)$ is a non-zero polynomial and $g(z)$ is a transcendental entire function such that $\sigma_{p}(w(z))=\sigma_p(e^{g(z)})=\infty$
and $\sigma_{p+1}(w(z))=\sigma_p(g(z))<\infty$.
If $w(z+1)-w(z-1)\equiv0$, then (\ref{5.1}) and (\ref{5.2}) give
\begin{equation*}
\frac{a_0(z)}{a(z)}-\frac{H'(z)}{H(z)}=g'(z),
\end{equation*}
a contradiction. So $w(z+1)-w(z-1)\not\equiv0$. Substituting (\ref{5.2}) into (\ref{5.1}) and letting
\begin{equation}\label{5.3}
F=a_0(z)-a(z)\left(\frac{H'(z)}{H(z)}+g'(z)\right),
\end{equation}
we get
\begin{equation}\label{5.4}
F=H(z+1)e^{g(z+1)}-H(z-1)e^{g(z-1)}
\end{equation}
and
\begin{equation}\label{5.5}
F'=(H'(z+1)+H(z+1)g'(z+1))e^{g(z+1)}-(H'(z-1)+H(z-1)g'(z-1))e^{g(z-1)}.
\end{equation}

If $H(z+1)(H'(z-1)+H(z-1)g'(z-1))\not\equiv H(z-1)(H'(z+1)+H(z+1)g'(z+1))$, then
by (\ref{5.4}) and (\ref{5.5}), we get
\begin{align}\label{5.6}
&e^{g(z-1)}=\nonumber\\
&\frac{(H'(z+1)+H(z+1)g'(z+1))F-H(z+1)F'}{H(z+1)(H'(z-1)+H(z-1)g'(z-1))-H(z-1)(H'(z+1)+H(z+1)g'(z+1))}.
\end{align}
By Lemma 5.1, we have $\sigma_p(e^{g(z-1)})=\sigma_p(e^{g(z)})=\infty$ and the iterated p-order of the right hand side of (\ref{5.6})
is no more than $\sigma_p(g(z))<\infty$. This is a contradiction.

If $H(z+1)(H'(z-1)+H(z-1)g'(z-1))\equiv H(z-1)(H'(z+1)+H(z+1)g'(z+1))$, then by (\ref{5.5}) we get
\begin{equation}\label{5.7}
H(z+1)F'=(H'(z+1)+H(z+1)g'(z+1))(H(z+1)e^{g(z+1)}-H(z-1)e^{g(z-1)}).
\end{equation}
Obviously, $H'(z+1)+H(z+1)g'(z+1)\not\equiv0$.
By (\ref{5.4}) and (\ref{5.7}), we get
\begin{equation*}
\frac{F'}{F}=\frac{H'(z+1)}{H(z+1)}+g'(z+1).
\end{equation*}
So $F=CH(z+1)e^{g(z+1)}$, where $C$ is a non-zero constant. By (\ref{5.3}) we get
\begin{equation}\label{5.8}
CH(z+1)e^{g(z+1)}=a_0(z)-a(z)\left(\frac{H'(z)}{H(z)}+g'(z)\right).
\end{equation}
Comparing the iterated p-order of both sides of (\ref{5.8}), we get a contradiction.

So we
proved that the equation (\ref{1.5}) has no entire solutions with finite iterated order.

\section*{ Acknowledgements}

The first author is partly supported by Guangdong National Natural Science Foundation of China (No. 2016A030313745)
and Training Plan Fund of Outstanding Young Teachers of Higher Learning Institutions of Guangdong Province of China (No. Yq20145084602).
The second author  is supported by Guangdong National Natural Science Foundation of China (No.2014A030313422).

\end{document}